\documentclass{amsart}
\usepackage{amstext}
\usepackage{bbding}
\usepackage{pifont}
\usepackage{amssymb}
\newtheorem{thm}{Theorem}[section]
\newtheorem{cor}[thm]{Corollary}

\newtheorem{pro}[thm]{Proposition}
\newtheorem{defin}[thm]{Definition}

\begin{document}

\title{Zero entropy invariant measures for some skew product diffeomorphisms}
\author{Peng Sun}
\date{\today}
\maketitle

\begin{abstract}
In this paper we study some skew product diffeomorphisms with nonuniformly
hyperbolic structure along fibers. We show that there is an invariant measure
with zero entropy which has atomic conditional measures along fibers.
This gives affirmative answer for these diffeomorphisms to the
question suggested by Herman that
a smooth diffeomorphism of positive topological entropy fails to be uniquely
ergodic. The proof is based on some techniques analogous to those
developed by Pesin (\cite{Pe1}) and Katok (\cite{Ka1}, \cite{KM1}) with investigation
on some combinatorial properties of the projected return map on the base.

\end{abstract}

\tableofcontents

\section{Introduction}

Let $f$ be a $\mathrm{C}^{1+\alpha}$ ($\alpha>0$) diffeomorphism of a compact
 $s$-dimensional smooth manifold $M$ and $\mathrm{d}f:\mathrm{T}M\to\mathrm{T}M$
the derivative of $f$. $f$ preserves a Borel probability measure
$\mu$. For every $x$ in a set $\Lambda$ of full measure, the Lyapunov
exponent
$$\chi(v,f)=\lim_{n\to\infty}\frac{\ln\|\mathrm{d}f^nv\|}{n}$$
exists for every nonzero vector $v\in \mathrm{T}_xM$. This functional takes
on at most $s$ values on $\mathrm{T}_xM$  and is independent of $x\in\Lambda
$ if $\mu$ is ergodic.
If all Lyapunov exponents are nonzero, then $\mu$ is called a hyperbolic
measure. Smooth systems with hyperbolic measures are called nonuniformly
hyperbolic. The theory for studying such systems was developed by Pesin and
then combined with some powerful techniques by A. Katok to look for
invariant orbits and produced a number of profound results. These techniques
serve as cornerstones for our discussion. For all necessary definitions,
theorems and background facts relevant to this paper, one may see \cite{BP1}
for quick
reference or \cite{BP2} for detailed proofs.

In \cite{Ka1} Katok showed:

\begin{thm}\label{TKa1}
Let $f$ be a $\mathrm{C}^{1+\alpha}(\alpha>0)$ diffeomorphism of a compact
manifold M, and $\mu$ a Borel probability $f$-invariant hyperbolic measure.
Then
$$\overline{Per}(f)\supset\mathrm{supp}(\mu)$$and
$$\max(0, \limsup_{n\to\infty}\frac{\ln
P_n(f)}{n})\geq h_\mu(f)$$ 
Where $Per(f)$ is the set of all periodic points of $f$ and $P_n(f)$ the
number of periodic points of $f$ with period $n$. $h_\mu(f)$ is the metric
entropy with respect to $\mu$.
\end{thm}

In particular, if the manifold $M$ is 2-dimensional, then by Ruelle inequality
\cite{Ruelle}, every ergodic invariant measure $\mu$ with positive metric entropy must be hyperbolic. Taking also the variational principle into account, we
have

\begin{cor}\label{CKa1}
        For any $\mathrm{C}^{1+\alpha}(\alpha>0)$ diffeomorphism $f$ of a
        2-dimensional compact manifold with positive topological entropy,
        
        \begin{equation}\label{eq1}\limsup_{n\to\infty}\frac{\ln P_n(f)}{n}\geq h(f)\end{equation}
        
        Hence $f$ is not minimal or uniquely ergodic.
\end{cor}

In general, Equation (\ref{eq1}) is not true for high dimensional cases.
There can be no periodic orbit for a diffeomorphism with positive topological
entropy. Herman \cite{He1} constructed a remarkable example as following:

Consider the $\mathrm{C}^\infty$ map $A:\mathbb{T}^1\to\mathrm{SL}(2,\mathbb{R})$
defined by    \[A(\theta)=A_\theta=\left(\begin{array}{r r}
    \cos 2\pi\theta & -\sin 2\pi\theta \\ \sin 2\pi\theta & \cos
    2\pi\theta
    \end{array} \right)
    \left(\begin{array}{r r}\lambda & 0 \cr 0 &
    1/\lambda\end{array}\right)\]
    where $\lambda>1$ is a fixed number. Let $R_\alpha:\mathbb{T}^1\to\mathbb{T}^1$ be the rotation by
    $\alpha\in\mathbb{T}^1-(\mathbb{Q}/\mathbb{Z})$.

\begin{thm}\label{THe1}\textnormal{(Herman, \cite{He1})}
There is a dense $\mathrm{G}_\delta$ subset W of $\mathbb{T}^1$, such that
for every $\alpha\in W$, the smooth diffeomorphism $F_\alpha=(R_\alpha, A(\theta))$
on $\mathbb{T}^1\times\mathrm{SL}(2,\mathbb{R})/\Gamma$, given by $(\theta,
y)\mapsto(\theta+\alpha, A(\theta)\cdot y)$, is minimal and has positive
topological entropy.
\end{thm}

Herman's example prompted a fruitful research, for example, on generic linear
cocycles over compact systems.
The phenomenon he discovered
 turned out to be common for $\mathrm{SL}(2,\mathbb{R})$ extension
  over rotations \cite{AK}.  

However, the diffeomorphisms in Herman's example fail to be uniquely ergodic.
 We can find a measurable transformation $S:\mathbb{T}^1\to\mathrm{SL}(2,\mathbb{R})$
 such that for almost every $\theta\in\mathbb{T}^1$, $H_\theta=S_{\theta+\alpha}A_\theta S_\theta^{-1}=\left(\begin{array}{l l}
l_\theta & 0 \\ 0 & l_\theta^{-1}\end{array}\right)$
    is diagonal. Then for every measure $\tau$ preserved by the geodesic
    flow which corresponds to the left action by
    $G_t=\left(\begin{array}{c c} \exp(t/2) & 0 \\ 0 &
    \exp(-t/2)\end{array}\right)$,
    $\mu_\tau=\int\tau\circ S_\theta dm$ is $F_\alpha$-invariant, where $m$
    is the Lebesgue measure on $\mathbb{T}^1$. In particular, if $\tau$ is supported on a periodic
    orbit of the geodesic flow, then $h_{\mu_\tau}(F_\alpha)=0$.
    
Whether a smooth diffeomorphism of positive topological entropy can be uniquely
ergodic is still in question (For homeomorphisms, the answer is yes. See
for examples, \cite{BC}). We studied some
skew product diffeomorphisms and found some invariant measures similar to
those in Herman's example.

Let $(X,m)$ be a probability measure space. $g:X\to X$ is an invertible transformation
(mod 0) preserving $m$. $M$ is an
$l$-dimensional compact Riemannian manifold. For every $x\in X$, $h_x:M\to
M$ is
a $\mathrm{C}^{1+\alpha}$ diffeomorphism. Assume $f=(g,h_x)$ on $X\times
M$ preserves a measure $\mu=\int\nu_xdm$. Let $\mathrm{T}_yM=\mathrm{T}_p(\{x\}\times
M)$ for $y=(x,p)\in X\times M$. In this paper
We prove:

\begin{thm}\label{Tmain}\textnormal{(Main Theorem)}
If for almost every $y\in Y$, The Lyapunov exponent $\chi(v,f)\neq
0$ for all $v\in\mathrm{T}_yM\backslash\{0\}$, then f has an invariant measure whose conditional measure on each
 fiber is atomic.
\end{thm}

Now suppose that we have a $\mathrm{C}^{1+\alpha}$ diffeomorphism $f=(g,
h_x)$ on $X\times M$. Assume
that $h(g)=0$ and $M$ is 2-dimensional. If $h(f)>0$, we must have $h_\mu(f)>0$
for some ergodic invariant measure $\mu=\int\nu_xdm$. Then by
Ledrappier-Young's formula \cite{LY1}, the Lyapunov exponents along fiber
direction must be nonzero
almost everywhere. Hence by Theorem \ref{Tmain} $f$ has an
invariant measure with atomic conditional measures along fibers. The following
statement avoids any mention of exponents.

\begin{cor}\label{Cmain}
If $f$ has positive topological entropy, $g$ has zero topological entropy,
and $M$ is 2-dimensional, then $f$ has a measure of zero entropy and is not
uniquely ergodic.
\end{cor}

\section{Shadowing lemma}

Now that we have a $\mathrm{C}^{1+\alpha}$ diffeomorphism $f=(g,h_x)$ on $Y=X\times M$ that has nonzero exponents along fiber direction. We may assume that
$g$ is ergodic by considering an ergodic component. Almost all results in \cite{Ka1, KM1, Pe1} can be adapted in this setting with careful modification.
By considering the derivative $\mathrm{d}_yh_x=\mathrm{d}_ph_x$ for $y=(x,p)$ as a linear cocycle over $f$, we have:

\begin{thm}\label{Tregular}
        Assume $\dim M=l$. Denote by $B^{k}(r)$ the standard Euclidean r-ball
        in $\mathbb{R}^k$ centered at the origin. There exists a set $\Lambda_0\subset Y$ of full measure such that
        for every sufficiently small $\epsilon>0$ and some $\chi>0$:
   
   \begin{enumerate}
   
   \item There exists a tempered function $q:\Lambda_0\to(0,1]$ and a collection
   of embeddings $\Psi_y:B^l(q(y))\to\{x\}\times M$ for each $y=(x,p)\in\Lambda_0$
   such that $\Psi_y(0)=y$ and $e^{-\epsilon}<q(y)/q(f(y))<e^\epsilon$.
   
   \item There exist a constant $K>0$ and a measurable function $C:\Lambda_0\to\mathbb{R}$
   such that for $z_1,z_2\in B^l(q(y))$,
   $$K^{-1}d(\Psi_y(z_1),\Psi_y(z_2))\leq\|z_1-z_2\|\leq C(y)d(\Psi(z_1),\Psi(z_2))$$
   with $e^{-\epsilon}<C(f(y))/C(y)<e^\epsilon$.
   
   \item The map $f_y:=\Psi^{-1}_{f(y)}\circ f\circ\Psi_y:B^s(q(y))\times
   B^{l-s}(q(y))\to\mathbb{R}^l=\mathbb{R}^k\times\mathbb{R}^{l-k}$
   has the form $$f_y(u,v)=(A_yu+\eta_{2,y}(u,v), B_yv+\eta_{1,y}(u,v))$$
   where $\eta_{1,y}(0,0)=\eta_{2,y}(0,0)=0$, $\mathrm{d}\eta_{1,y}(0,0)=\mathrm{d}\eta_{2,y}(0,0)=0$
   and
   $$\|A_y\|<\exp-(\chi-\epsilon),\|B_y^{-1}\|<\exp-(\chi-\epsilon)$$
   For $z=(u,v)\in B^l(q(y))$, $\eta_y(z)=(\eta_{1,y}(z),\eta_{2,y}(z))$:
   $$\|\mathrm{d}_z\eta_y\|<\epsilon, \|\eta_y(z)\|<\epsilon$$
   \end{enumerate}
\end{thm}

\begin{defin}\label{Dregular}
The points $y\in\Lambda_0$ are called regular points. For each regular point
$y$, the set $N(y)=\Psi_y(B(q(y)))$
is called a regular neighborhood of $y$. Let $r(y)$ be the radius of the
maximal ball contained in the regular neighbor hood $N(y)$. We say $r(y)$
is the size of $N(y)$.
\end{defin}

\begin{thm}\label{Tdelta}
For each $\delta>0$ and each sufficiently small $\epsilon(\delta)>0$, there is a
set $\Lambda_\delta\subset\Lambda_0$ which has compact intersection $\Lambda_{\delta,x}$
(may be empty)
with each fiber $\{x\}\times M$, such that $\mu(\Lambda_\delta)>1-\delta$
and the following conditions hold:

\begin{enumerate}
\item The functions $y\mapsto q(y)$, $y\mapsto C(y)$ and $y\mapsto\Psi_y$
as in Theorem \ref{Tregular} for $\epsilon=\epsilon(\delta)$,
and $y\mapsto r(y)$ are all continuous on $\Lambda_{\delta,x}$ for each $x\in
X$.

\item The decomposition $\mathrm{T}_yM=\mathrm{d}_y\Psi_y\mathbb{R}^k\times\mathrm{d}_y\Psi_y\mathbb{R}^{l-k}$
depends continuously on $y$ in $\Lambda_{\delta,x}$.

\item On $\Lambda_\delta$, there are bounds: $q_\delta=\min\{q(y)\}$, $r_\delta=\min\{r(y)\}$,
$C_\delta=\max\{C(y)\}$.

\end{enumerate}
\end{thm}

With similar definitions and properties for admissible manifolds, we are
able to derive the following version of Shadowing Lemma:

\begin{thm}\label{Tshadowing}
Given $\delta>0$, for $\bar q<q_\delta$, set $\tilde{\Lambda}_\delta(\bar
q)=\bigcup_{y\in\Lambda_\delta}\Psi_y(B(0,\bar q))$. Given $a\in\mathbb{Z}\bigcup\{-\infty\}$
and $b\in\mathbb{Z}\bigcup\{\infty\}$, a sequence $\{y_n=(x_n,p_n)\}_{a<n<b}$ is called an
$(\delta,\bar q)$-pseudo orbit for $f=(g,h_x)$ if there are $\{z_n\in\Lambda_{\delta,x_n}\}_{a<n<b}$
and $\{k_n\}_{a<n<b}$ such that
for every $n$, $y_{n}\in\Psi_{z_n}(B(0,\bar q))$ and $f^{k_{n+1}-k_n}(y_n)\in\Psi_{z_{n+1}}(B(0,\bar q))$.
Then there exists $\gamma=\gamma(\delta)$
such that for every $(\delta,\gamma)$-pseudo orbit, there is a unique point
$\tilde y\in Y$ such that $f^{k_n}(\tilde y)\in\Psi_{z_n}(B(0,q_\delta))$ for
all $a<n<b$.
\end{thm}
 
\section{Integrability of return time}

Now we would like to take a proper Pesin set on which the shadowing techniques
can be carried out.

\begin{defin}\label{defTube}
Let $\pi:Y\to X$ be the projection to the base.
A measurable subset $P\subset Y$ is called a "Regular Tube", if for some $\delta>0$, $\epsilon>0$, $\nu_0>0$ and $\gamma=\gamma(\delta)$
as in Theorem \ref{Tshadowing}, there exists for every $x\in B=\pi(P)$, a
point $z(x)\in\Lambda_{\delta,x}$ such that $P_x=P\bigcap(\{x\}\times M)\subset\Psi_{z(x)}(B(0,\gamma))$,
$\nu_x(P_x)>\nu_0$ and $m(B)>1-\epsilon$.
\end{defin}

The existence of such a "Regular Tube" is guaranteed in Section 2. In this
"Regular Tube", we
can take a measurable section $s:B\to P$, $\pi\circ s=\mathrm{id}_B$. Let $S=s(B)$. We are then going
to consider the first return map $f_P$ on $P$.

\begin{pro}\label{Treturn}
$s$ can be chosen in such a way that the first return time from $S$ to $P$
is integrable with respect to $m$. In particular, we may assume every point
in $S$ returns to $P$ in finite times.
\end{pro}

{\flushleft\emph{Proof.}} For every $y\in P$, denote by $n(y)$ the return time
of $y$. Since $\mu$ is $f$-invariant and $\mu(P)>m(B)\cdot\nu_0>0$, we have:
$$0<\int_{P}n(y)\mathrm{d}\mu=\mu(\bigcup_{j\ge0}F^j(P'))\le1$$
But
$$\int_{P}n(y)\mathrm{d}\mu=\sum^\infty_{j=0}\mu(P_j)$$
where $P_j=P\backslash(\bigcup_{1\le k\le j}F^{-k}(P))$.

We may choose $s$ such that for every $x\in B$, $s(x)\in P_j$ only if $\nu_x(P_x\backslash P_j)=0$. Let $B_j=\pi_1(S\backslash(\bigcup_{1\le k\le j}F^{-k}(P))=\pi_1(S\bigcap P_j)$. By the way $s$ is chosen and the assumption $\nu_x(P_x)>\nu_0$, we
have $\mu(P_j)\geq\mu(\pi_1^{-1}(B_j)\bigcap P)>m(B_j)\cdot\nu_0$, hence
$$\int_{B}n(s(x))\mathrm{d}m=\sum^\infty_{j=0}m(B_j)<\sum^\infty_{j=0}\frac{1}{\nu_0}\mu(P_j)\le\frac{1}{\nu_0}$$
The return time is integrable. \SquareSolid

\section{Projected return map on the base}

Now let us try to find the invariant measure described in Theorem \ref{Tmain}.
We may assume that $g$ has no periodic point,
or else the problem can be reduced to the case considered by Katok \cite{Ka1}.
Moreover, $g$ is invertible as we assumed earlier.

We are looking for an invariant set $I$ which has finite intersection $I_x$ with almost every fiber $\{x\}\times M$. The measure $\tau_x$ supported on $I_x$ is the delta counting
measure on $\{x\}\times M$. Then an invariant measure for $f$ can be given by $\int\tau_xdm$.

To find the invariant set $I$, we start with a fixed "Regular Tube" $P$ and a measurable
section $s$ as specified in section 3. 
Let $f_P$ be the first return map on $P$ and $g_B$
the first return map for $g$ on $B$. Define
$r:B\to B$ by $r(x)=\pi\circ f_P(s(x))$. $r(x)$ is the projection of the
return map on the base. $g_B$ is invertible but $r$ may
not. For every $x\in B$, let $k(x)$ be such that $f^{k(x)}=f_P(s(x))$.

We can define a partial order on $B$: $x_1\prec x_2$ iff there
is $n\geq 0$ such that $g_B^n(x_1)=x_2$, i.e. $x_2$ is an image of $x_1$ under
iterates of $g$. Since $g$ is invertible and has no periodic point, this partial
order is well defined. Moreover, if there is $n\geq 0$ such that $r^{n}(x_1)=x_2$
then we write $x_1\prec\prec x_2$, which implies $x_1\prec x_2$. This is
also a partial order.

We can define an equivalence relation on $B$: $x_1\backsim x_2$ iff $Q(x_1,x_2):=\{x\in
B|x_1\prec\prec x\text{ and }x_2\prec\prec x\}\neq\emptyset$, i.e.
there are $n_1, n_2>0$ such that $r^{n_1}(x_1)=r^{n_2}(x_2)$. If $x\in Q(x_1,x_2)$
then $r^n(x)\in Q(x_1,x_2)$ for all $n>0$. If $x_1\backsim
x_2$, we must have  $x_1\prec x_2$ or $x_2\prec x_1$, denoted by $x_1\precsim
x_2$ or $x_2\precsim x_1$. If $x_1\precsim x_2$, define
$\sigma(x_1, x_2)$ as the minimal (with respect to $\prec$, always in this
paper) element in $Q(x_1,x_2)$. In particular,
if $x_1\prec\prec
x_2$ then $x_1\precsim x_2$ and $\sigma(x_1,x_2)=x_2$.

{\flushleft\emph{Remark.}} The equivalence relation $\sim$ defined here is
crucial in this paper.
If $x_1\precsim x_2$, then $s(x_1)$ and $s(x_2)$ return
to the same fiber after iteration of $f_P$. However, $s(x_1)$ does not necessarily
return to $P_{x_2}$, i.e. two points in $S$ may return to $P$ on the same
fiber and this may happen all the time. We had trouble dealing with this situation while looking for pseudo orbits. Introduction of this equivalence relation solved
this problem. We can then, in each equivalence class, find a unique
orbit of $r$
(a sequence of returns, lifted to a pseudo orbit) to construct the invariant set.

\begin{pro}\label{Pro2}
For almost every $x\in B$, $J(x)=\{\sigma(x',x)|x'\precsim x,x'\neq x\}$
is finite. Denote
by $x^*$ the maximal element of $J(x)$. Then $x'\prec\prec x^*$ for all
$x'\precsim x$. Let $W(x)=\{\bar x|x'\prec\prec\bar x\text{ for all }x'\precsim
x\}$, then $x^*=\min W(x)$. Moreover, if $x_1\precsim x_2$, then $x_1^*\prec\prec x_2^*$.
\end{pro}

{\flushleft\emph{Proof.}} For every $x\in B$, define the set of "jumps" $J'(x):=\{r(x')|x'\prec x, x'\ne x\text{ and }x\prec
r(x')\}$. By integrability of return times (Proposition \ref{Treturn}), $J'(x)$ must be finite for almost
every $x\in B$. To see this, we can consider the set $$\tilde S=\bigcup_{x\in B}\{f(s(x)),
f^2(s(x)),\cdots,f^{k(x)}(s(x))\}$$
Let $\tilde S_j=\{x\in X||\tilde S\bigcap(\{x\}\times
M)|=j\}$. We can count the return times and
get
$$\sum^\infty_{j=0}j\cdot m(\tilde S_j)=\int_Bk(x)\mathrm{d}m<\infty$$
and $|J'(x)|\le|\tilde S\bigcap(\{x\}\times M)|<\infty$ for almost every
 $x\in B$.

For every $x'\precsim x$ but $x'\ne x$, there must be $\bar x\in J'(x)$ such that
$x'\prec\prec\bar x$ and $\sigma(x',x)=\sigma(x,\bar x)$. So for different
elements $x_1,x_2\in J(x)$, there must be different elements $\bar x_1, \bar x_2\in J'(x)$ such that $\sigma(\bar x_i,x)=x_i$, $i=1,2$. Hence $|J(x)|\leq|J'(x)|<\infty$.

By definition, $x'\prec\prec\sigma(x',x)$ for every $x'\in B$, and
$\sigma(x_1,x)\prec\prec\sigma(x_2,x)$ if $\sigma(x_1,x)\prec\sigma(x_2,x)$
since they are both images of $x$ under iteration of $r$. So for every $x'\precsim
x$, we have $x'\prec\prec\sigma(x',x)\prec\prec x^*$.

Since $x^*\in J(x)$, there is some $x'\precsim x$ such that $\sigma(x',x)=x^*$.
Then for every $\bar x\in W(x)$, $\bar x\in Q(x',x)$. But $x^*=\sigma(x',x)=\min
Q(x',x)\prec\bar x$. $x^*=\min W(x)$.

If $x_1\precsim x_2$, then $x_1^*=\min\{x|x'\prec\prec x\text{ for all }x'\precsim
x_1\}\prec\min\{x|x'\prec\prec x\text{ for all }x'\precsim x_2\}=x_2^*$,
because  the second set is contained in the first one. But $x_1\prec\prec
x_1^*$ and $x_1\prec\prec
x_2^*$, from previous discussion we must have $x_1^*\prec\prec x_2^*$.
\SquareSolid

\begin{pro}\label{Pro3}
 Let $B_0=\{x\in B|\text{there is no such }x'\neq x\in B\text{ that }x'\precsim
 x\}$. Then $m(B_0)=0$. Hence by replacing $B$ by $B\backslash(\bigcup_{k\in\mathbb{Z}}g_B^k(B_0))$
and $P$ accordingly, we may assume that for every $x\in B$, there is at least
one element $x'\in B$ such that
$x'\precsim x$ but $x'\neq x$.
\end{pro}

{\flushleft\emph{Proof.}} If $m(B_0)>0$, then there must be an element $x_0\in B_0$ such that
$B_0(x_0)=\{g_B^{-n}(x_0), n\in\mathbb{N}\}\bigcap B_0$ has infinitely many elements
by Poincar\'e Recurrence Theorem because $g_B$ is invertible and $m$-preserving.
From the proof of Proposition \ref{Pro2}, $J'(x_0)$ has finitely many elements.
But for every $x\in B_0(x_0)$, there must be $x'\in J'(x_0)$ such that $x\prec\prec
x'$. Hence there is an element $\tilde x\in J'(x_0)$ such that $\tilde B_0(x_0)=\{x\in
B_0(x_0)|x\prec\prec\tilde x\}$ has infinitely many elements. But $x_1\sim
x_2$ for all $x_1,x_2\in\tilde B_0(x_0)\subset B_0$ because $\bar x\in Q(x_1,x_2)\neq\emptyset$, which is a contradiction.
\SquareSolid

For every $x\in B$, let $G(x):=\{x'\in B|x'\sim x\}$ and $G^*(x):=\{(x')^*|x'\in
G(x)\}$. If $x_1\sim x_2$, then we must have $G(x_1)=G
(x_2)$ and $G^*(x_1)=G^*(x_2)$.

 Pick $H(x)$ in the following way:

\begin{enumerate}
\item If $G^*(x)$ is not properly defined: Let $H(x)=\emptyset$. (only for
$x$ in a set of measure zero)
\item
If $G^*(x)$ has a minimal element $\tilde x$: Let $H(x)=\{f^n(s(\tilde
x))\}_{-\infty<n<\infty}$.
\item
If $G^*(x)$ has no minimal element: By Proposition \ref{Pro2}, $G^*(x)$ can
be completed to a full orbit
of $r$, $\bar G(x)=\bigcup_{0\le n<\infty}r^n(G^*(x))$. In each equivalence class $G(x)$, $\bar G(x)$ is a sequence of
returns and is uniquely defined in the sense $\bar G(x)=\bigcup_{x_1\sim
x}\bigcap_{x_2\precsim x_1}\{r^n(x_2)\}_{0\le n<\infty}$. $\bar G(x)$ ordered
by "$\prec\prec$" can be viewed as a sequence $\{\bar x_n\}_{-\infty<n<\infty}$
and $r(x_n)=x_{n+1}$ for all $n$. Then the sequence $\{s(\bar x_n)\}_{-\infty<n<\infty}$
is in fact a $(\delta,\gamma)$-pseudo orbit. Let us call it the pseudo orbit
associated to $x$. Note the pseudo orbits associated to equivalent elements
coincide.
 By the way the "Regular Tube" $P$
was chosen, we can find $\tilde y\in Y$ as specified in Theorem \ref{Tshadowing}.
Let $H(x)=\{f^n(\tilde y)\}_{-\infty<n<\infty}$.
\end{enumerate}

Let $I=\bigcup_{x\in B}H(x)$. By definition, $H(x)$ is invariant
for all $x\in B$. Hence $I$ is $f$-invariant.

\begin{pro}\label{Profinite}
For almost every $x\in X$, $I_x=I\bigcap(\{x\}\times M)$ is nonempty and
contains finitely many elements.
\end{pro}

{\flushleft\emph{Proof.}} 
For almost every $x\in B$, $I_x\supset(H(x)\bigcap(\{x\}\times M))\ne\emptyset$
by definition.
Note that $H(x_1)=H(x_2)$ if $x_1\sim x_2$.
For different elements $y_1,y_2\in I_x$, there are $x_1, x_2\in B$ such that
$y_i\in H(x_i)$, $i=1,2$. We must have $x_1\nsim x_2$ and $G(x_1)\bigcap
G(x_2)=\emptyset$. But $G(x_i)\bigcap J'(x)\ne\emptyset$, $i=1,2$.
So we have $|I_x|\leq|J'(x)|$.

Recall that $m(B)>0$, so by ergodicity, $I_x$ must be nonempty
and contain
finitely many elements for almost every $x\in X$. \SquareSolid

Since $I$ is invariant and  $I_x$ is finite, $\int\tau_x dm$ is an invariant
measure as requested, where $\tau_x$ is the delta counting measure on $\{x\}\times
M$ supported on $I_x$, for almost every $x\in X$. The entropy of this measure
is the same as the entropy of the transformation $g$ on the base.

\section{Measures of intermediate entropies}

In \cite{Ka2}, Katok showed a stronger result:

\begin{thm}\label{K2}
If $f:M\to M$ is a $\mathrm{C}^{1+\alpha}$ diffeomorphism of a compact smooth
manifold and $\mu$ an ergodic hyperbolic measure for $f$ with $h_\mu(f)>0$,
then for any $\epsilon>0$ there exists a hyperbolic horseshoe $\Gamma$ such
that $h(f|_\Gamma)>h_\mu(f)-\epsilon$. Hence for any number $\beta$ between
zero and $h_\mu(f)$,
there is an ergodic invariant measure $\mu_\beta$ such that $h_{\mu_\beta(f)}=\beta$.
\end{thm}

We are looking for analogous result to the theorem for our skew product diffeomorphisms.
 As we know, for skew
product diffeomorphisms there may not be any proper closed invariant sets.
 However,
 we may expect
to have an invariant set $\Gamma$ that has closed intersection with almost
every fiber, on which $f$ acts like a horseshoe map.
Moreover, this horseshoe should carry
an entropy arbitrarily close to $h_\mu(f)$ in order to produce invariant
measures with arbitrary intermediate entropies. This work is in progress.

We may also ask the question if theorem  \ref{K2} holds for any $\mathrm{C}^{1+\alpha}$
 diffeomorphism without the assumption that $\mu$ is hyperbolic. 
 We have not yet found even
 a zero entropy measure in this general case.



\end{document}